\documentclass[11pt]{article}
\usepackage[a4paper]{anysize}\marginsize{1.5cm}{1.5cm}{1.5cm}{1cm}
\pdfpagewidth=\paperwidth \pdfpageheight=\paperheight
\usepackage{amsfonts,amssymb,amsthm,amsmath,eucal,tabu,url}
\usepackage{pgf}
 \usepackage{array}
 \usepackage{pstricks}
 \usepackage{pstricks-add}
 \usepackage{pgf,tikz}
 \usetikzlibrary{automata}
 \usetikzlibrary{arrows}
 \usepackage{indentfirst}
\theoremstyle{plain}
\newtheorem{thm}{Theorem}[section]
\newtheorem{theorem}[thm]{Theorem}
\newtheorem*{theoremA}{Theorem A}
\newtheorem*{theoremB}{Theorem B}
\newtheorem*{theoremC}{Theorem C}
\newtheorem*{theoremD}{Theorem D}

\theoremstyle{definition}
\newtheorem{definition}[thm]{Definition}
\newtheorem{remark}[thm]{Remark}
\newtheorem{example}[thm]{Example}

\newtheorem{question}[thm]{Question}

\newtheorem{thevarthm}[thm]{\varthmname}

\newenvironment{varthm*}[1]{\trivlist\item[]{\bf #1.}\it}{\endtrivlist}

\renewcommand\ge{\geqslant}
\renewcommand\geq{\geqslant}

\renewcommand\leq{\leqslant}

\newcommand\be{\begin{eqnarray*}}
\newcommand\ee{\end{eqnarray*}}

\newcommand\newop[2]{\def#1{\mathop{\rm #2}\nolimits}}
\newop\log{log}
\newop\ord{ord}
\newop\Gal{Gal}
\newop\SL{SL}
\newop\Bl{Bl}
\newop\mult{mult}
\newop\mass{mass}
\newop\div{div}
\newop\codim{codim}
\newop\sing{sing}
\newop\vdim{vdim}
\newop\edim{edim}
\newop\Ass{Ass}
\newop\size{size}
\newop\reg{reg}
\newop\satdeg{satdeg}
\newop\supp{supp}
\newop\Neg{Neg}
\newop\Nef{Nef}
\newop\Nefh{Nef_H}
\newop\Eff{Eff}
\newop\Zar{Zar}
\newop\MB{MB}
\newop\MBxC{MB\mathit{(x,C)}}
\newop\NnB{NnB}
\newop\Bigg{Big}
\newop\Effbar{\overline{\Eff}}

\def\keywordname{{\bfseries Keywords}}%
\def\keywords#1{\par\addvspace\medskipamount{\rightskip=0pt plus1cm
\def\and{\ifhmode\unskip\nobreak\fi\ $\cdot$
}\noindent\keywordname\enspace\ignorespaces#1\par}}
\def\subclassname{{\bfseries Mathematics Subject Classification
(2000)}\enspace}
\def\subclass#1{\par\addvspace\medskipamount{\rightskip=0pt plus1cm
\def\and{\ifhmode\unskip\nobreak\fi\ $\cdot$
}\noindent\subclassname\ignorespaces#1\par}}

\begin{document}
\title{On line and pseudoline configurations and ball-quotients}
\author{J\"urgen Bokowski, Piotr Pokora}

\date{\today}
\maketitle
\thispagestyle{empty}
\begin{abstract}
In this note we show that there are no real configurations of $d\geq 4$ lines in the projective plane such that the associated Kummer covers of order $3^{d-1}$ are ball-quotients and there are no configurations of $d\geq 4$ lines such that the Kummer covers of order $4^{d-1}$ are ball-quotients. Moreover, we show that there exists only one configuration of real lines such that the associated Kummer cover of order $5^{d-1}$ is a ball-quotient.  In the second part we consider the so-called topological $(n_{k})$-configurations and we show, using Shnurnikov's inequality, that for $n < 27$ there do not exist $(n_{5})$-configurations and and for $n < 41$ there do not exist $(n_{6})$-configurations.

\subclass{14C20, 52C35, 32S22}
\end{abstract}

%*****************************************************************************
\section{Preliminaries}
In his pioneering paper Hirzebruch \cite{Hirzebruch} constructed some new examples of algebraic surfaces which are ball-quotients, i.e., surfaces of general type satisfying equality in the Bogomolov-Miyaoka-Yau inequality \cite{M84}
$$K_{X}^{2} \leq 3e(X),$$
where $K_{X}$ denotes the canonical divisor and $e(X)$ is the topological Euler characteristic. The key idea of Hirzebruch, which enabled constructing these new ball-quotients, is that one can consider abelian covers of the complex projective plane branched along line configurations. Let us recall briefly how the celebrated construction of Hirzebruch works (for more details please consult for instance \cite{BHH87}).

Let $\mathcal{L} = \{l_{1}, ..., l_{d} \} \subset \mathbb{P}^{2}$ be a configuration of $d \geq 4$ lines such that there is no point $p$ where all $d$-lines meet and pick $n\in\mathbb{Z}_{\geq 2}$.
Now we can consider the Kummer extension having degree $n^{d-1}$ and Galois group $(\mathbb{Z}/n\mathbb{Z})^{d-1}$ defined as the function field
$$K: = \mathbb{C}\left(z_{1}/z_{0}, z_{2}/z_{0}\right)\left((l_{2}/l_{1})^{1/n}, ...,(l_{d}/l_{1})^{1/n}\right)$$
This Kummer extension is an abelian extension of the function field of the complex projective plane.
It can be shown that $K$ determines an algebraic surface $X_{n}$ with normal singularities which ramifies over the plane with the arrangement as the locus of the ramification. Hirzebruch showed that $X_{n}$ is singular exactly over a point $p$ iff $p$ is a point of multiplicity $\geq 3$ in $\mathcal{L}$. After blowing up these singular points we obtain a smooth surface $Y_{n}^{\mathcal{L}}$. It turns out that the Chern numbers of $Y_{n}^{\mathcal{L}}$ can be read off directly from combinatorics of line configurations, i.e.
$$\frac{c_{2}(Y_{n}^{\mathcal{L}})}{n^{d-3}} = n^{2}(3-2d+f_{1}-f_{0}) + 2n(d-f_{1}+f_{0}) + f_{1}-t_{2},$$
$$\frac{c_{1}^{2}(Y_{n}^{\mathcal{L}})}{n^{d-3}} = n^{2}(-5d+9+3f_{1}-4f_{0}) + 4n(d-f_{1}+f_{0}) +f_{1}-f_{0}+d+t_{2},$$
where $t_{r}$ denotes the number of $r$-fold points (i.e. points where exactly $r$ lines meet), $f_{0} = \sum_{r\geq 2}t_{r}$ and $f_{1} = \sum_{r\geq 2}r t_{r}$.
Moreover, it can be shown that $Y_{n}^{\mathcal{L}}$ has non-negative Kodaira dimension if $t_{d}=t_{d-1} = t_{d-2} = 0$ and $n \geq 2$, or $t_{d}=t_{d-1}=0$ and $n \geq 3$ (we assume additionally that $d \geq 6$), and in these cases we have $K_{Y_{n}^{\mathcal{L}}}^{2} \leq 3e(Y_{n}^{\mathcal{L}})$. Now we can define the following Hirzebruch polynomial (for more details, please consult the original paper due to Hirzebruch \cite[Section 3.1]{Hirzebruch}):
\begin{equation}
\label{Hirzpoly}
P_{\mathcal{L}}(n) = \frac{3e(Y_{n}^{\mathcal{L}}) - K_{Y_{n}^{\mathcal{L}}}^{2}}{n^{d-3}} = n^{2}(f_{0} - d) + 2n(d-f_{1}+f_{0}) + 2f_{1} + f_{0} - d - 4t_{2}
\end{equation}
and by the construction $P_{\mathcal{L}}(n) \geq 0$ provided that $n \geq 2$. If there exists a configuration of lines $\mathcal{A}$ such that there exists $m \in \mathbb{Z}_{\geq 2}$ with $P_{\mathcal{A}}(m)=0$, then $Y_{m}^{\mathcal{A}}$ is a ball quotient. There are some examples of line configurations which allow us to construct ball quotients via Hirzebruch's construction.
\begin{example}(\cite[p. 133]{Hirzebruch})
Let us consider the following configuration, which is denoted in the literature by $\mathcal{A}_{1}(6)$.
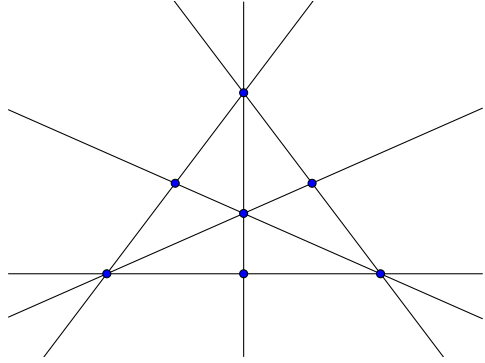
\begin{figure}[h]
\centering
\definecolor{qqqqff}{rgb}{0.,0.,1.}
\begin{tikzpicture}[line cap=round,line join=round,>=triangle 45,x=1.0cm,y=1.0cm,scale=0.6]
\clip(-5.155403622953033,-1.829197946106907) rectangle (5.234321773838792,6.006628470467939);
\draw [domain=-5.155403622953033:5.234321773838792] plot(\x,{(--12.--4.*\x)/3.});
\draw [domain=-5.155403622953033:5.234321773838792] plot(\x,{(--12.-4.*\x)/3.});
\draw [domain=-5.155403622953033:5.234321773838792] plot(\x,{(-0.-0.*\x)/6.});
\draw (0.,-1.829197946106907) -- (0.,6.006628470467939);
\draw [domain=-5.155403622953033:5.234321773838792] plot(\x,{(--6.--2.*\x)/4.5});
\draw [domain=-5.155403622953033:5.234321773838792] plot(\x,{(-6.--2.*\x)/-4.5});
\begin{scriptsize}
\draw [fill=qqqqff] (-3.,0.) circle (2.5pt);
\draw [fill=qqqqff] (3.,0.) circle (2.5pt);
\draw [fill=qqqqff] (0.,4.) circle (2.5pt);
\draw [fill=qqqqff] (-1.5,2.) circle (2.5pt);
\draw [fill=qqqqff] (1.5,2.) circle (2.5pt);
\draw [fill=qqqqff] (0.,0.) circle (2.5pt);
\draw [fill=qqqqff] (0.,1.3333333333333333) circle (2.5pt);
\end{scriptsize}
\end{tikzpicture}
\end{figure}
\noindent

Simple computations give 
$$P_{\mathcal{A}_{1}(6)}(n) = n^{2} - 10n +25,$$
which means that $Y_{5}^{\mathcal{A}_{1}(6)}$ is a ball-quotient.
\end{example}
\begin{example}(\cite[p. 133]{Hirzebruch})
Let us now consider the Hesse configuration $\mathcal{H}$ of lines (which cannot be drawn over the real numbers) having the following combinatorics:
$$d = 12, t_{2} = 12, t_{4} = 9.$$
Then
$$P_{\mathcal{H}}(n) = 9(n^{2} -6n + 9),$$
which means that $Y_{3}^{\mathcal{H}}$ is a ball-quotient. 
\end{example}
It is known that there are only a few examples of ball-quotients provided by line arrangements and it seems to be extremely difficult to find other examples. In this note we study a natural question about the existence of new ball quotients constructed via Hirzebruch's method. Before we formulate our main results let us define the following object.
\begin{definition}
Let $Y_{n}^{\mathcal{L}}$ be the minimal desingularization of $X_{n}$ constructed as the Kummer extension. Then $Y_{n}^{\mathcal{L}}$ is called the Kummer cover of order $n^{d-1}$.
\end{definition}
\begin{question}
Does a real line configuration $\mathcal{L} \subset \mathbb{P}_{\mathbb{C}}^{2}$ exist such that $Y_{3}^{\mathcal{L}}$ is a ball quotient?
\end{question}
\begin{remark}
In this note by a real line configuration we mean a configuration of lines which is realizable over the real numbers. For instance, the Hesse line configuration is not realizable over the real numbers.
\end{remark}
Our main results of this paper are the following strong classification results (our proofs are purely combinatorial).
\begin{theoremA}
There does not exist any real line configuration $\mathcal{L}$ with $d\geq 4$ lines and $t_{d}=t_{d-1}=0$ such that $Y_{3}^{\mathcal{L}}$ is a ball quotient.
\end{theoremA}

\begin{theoremB}
There does not exist any line configuration $\mathcal{L}$ with $d \geq 4$ lines and $t_{d} = t_{d-1} = 0$ such that $Y_{4}^{\mathcal{L}}$ is a ball-quotient.
\end{theoremB}
As a simple application of our methods we show the following results.
\begin{theoremC}
The configuration $\mathcal{A}_{1}(6)$ is (up to projective equivalence) the only configuration for $d\geq 4$ real lines such that the Kummer cover of order $5^{d-1}$ is a ball quotient.
\end{theoremC}

In our proof of Theorem A we use, in a very essential way, Shnurnikov's inequality (\ref{Shnurnikov}) for pseudoline configurations. Using this inequality we can prove the following result about topological $(n_{k})$-configurations.
\begin{theoremD}
For $n < 27$ there does not exist a topological $(n_{5})$-configuration and for $n < 41$ there does not exist a topological $(n_{6})$-configuration
\end{theoremD}
\section{Real line configurations and ball-quotients}

Firstly, we recall that the Hirzebruch polynomial, depending on $n \in \mathbb{Z}_{\geq 2}$, parameterizes the whole family of Hirzebruch's inequalities. Taking this into account, observe that if $n=3$, then we have the following inequality (we assume here that $t_{d} = t_{d-1} = 0$):
\begin{equation}
\label{weakHirz}
t_{2} + t_{3} \geq d + \sum_{r\geq 5} (r-4)t_{r}.
\end{equation}
It is worth pointing out that in a subsequent paper on the topic \cite{Hirzebruch1} Hirzebruch has improved his inequality (here we assume that $t_{d} = t_{d-1} = t_{d-2} = 0$):
\begin{equation}
\label{strongHirz}
t_{2} + \frac{3}{4}t_{3} \geq d + \sum_{r\geq 5}(2r-9)t_{r},
\end{equation}
and we should notice that this improvement comes from the Hirzebruch polynomial for $n=2$ with some extra effort -- please consult \cite{Hirzebruch1} for further details.

We will also need the following Melchior's inequality, which is true for real line configurations with $d\geq 3$ lines and $t_{d}=0$:
\begin{equation}
\label{Melchior}
t_{2} \geq 3 + \sum_{r\geq 4}(r-3)t_{r}.
\end{equation}

Finally, let us recall the notion of (real) pseudoline configurations.
\begin{definition}
We say that $\mathcal{C} \subset \mathbb{P}^{2}_{\mathbb{R}}$ is a configuration of pseudolines if it is a configuration of $n\geq 3$ smooth closed curves such that
\begin{itemize}
\item every pair of pseudolines meets exactly once at a single crossing (i.e., locally this intersection looks like $xy=0$),
\item curves do not intersect simultaneously at a single point.
\end{itemize}
\end{definition}
In particular, every real line configuration is a pseudoline configuration. Recently I. N. Shnurnikov \cite{Shnurnikov} has shown the following beautiful inequality.

\begin{theorem}
\label{Schnur}
Let $\mathcal{C}$ be a configuration of $n$ pseudolines such that $t_{n} = t_{n-1} = t_{n-2} = t_{n-3} = 0$. Then
\begin{equation}
\label{Shnurnikov}
t_{2} + \frac{3}{2} t_{3} \geq 8 + \sum_{r\geq 4} (2r-7.5)t_{r}.
\end{equation}
\end{theorem}

Now we are ready to prove Theorem A.
\begin{proof}
Our problem boils down to show that there does not exist a real line configuration satisfying
\begin{equation}
\label{equality}
t_{2} + t_{3} = d + \sum_{r\geq 5}(r-4)t_{r}.
\end{equation}
We start with excluding the case of $t_{d-2} = 1$ for which two possibilities remain (we assume here that $d\geq 6$)
\begin{itemize}
\item $\mathcal{A}_{1} : t_{d-2} = 1, t_{2} = 2d-3$,
\item $\mathcal{A}_{2} : t_{d-2} = 1, t_{3} = 1, t_{2} = 2d - 6$,
\end{itemize}
but it is easy to see that $\mathcal{A}_{1}$ and $\mathcal{A}_{2}$ do not satisfy (\ref{equality}).

From this point on we consider only real line configurations with $d$ lines where $t_{d} = t_{d-1} = t_{d-2} = 0$. Assume there exists a real line configuration $\mathcal{L}$ such that $Y_{3}^{\mathcal{L}}$ is a ball-quotient. Using (\ref{strongHirz}) and (\ref{equality}) we obtain
$$ -\frac{1}{4} t_{3} \geq \sum_{r\geq5} (r-5)t_{r},$$
which means that if $d\geq 4$ we have $t_{2} \geq 3$, $t_{3} = 0$ and $t_{r} = 0$ for $r\geq 6$. Moreover, it might happen that $t_{4}$ or $t_{5}$ are non-zero. This reduces (\ref{equality}) to
$$t_{2} = d + t_{5}.$$
On the other hand, we have the following combinatorial equality
$$d(d-1) = \sum_{r\geq 2} r(r-1) t_{r} = 2t_{2} + 12t_{4} + 20t_{5},$$
and combining this with $t_{2} = d + t_{5}$ we obtain
$$d(d-3) = 12t_{4} + 22t_{5}.$$
Using (\ref{Melchior}) we get
$$d-3 \geq t_{4} + t_{5}$$ and finally
$$12t_{4} + 22 t_{5} = d(d-3) \geq d(t_{4} + t_{5}),$$
which leads to
$$d \leq \frac{12t_{4} + 22t_{5}}{t_{4} + t_{5}} \leq 22.$$

Summing up, $\mathcal{L}$ satisfies the following conditions:
$$d \in \{4, ... ,22\}, \quad t_{2} = d + t_{5}, \quad d(d-3) = 12t_{4} + 22t_{4}, \quad d-3 \geq t_{4} + t_{5}.$$

It can be checked (for instance using a computer program) that the above constrains result in the following combinatorics (using the following convention in our listing : $\mathcal{L} = [d, t_{4}, t_{5}]$):
$$\mathcal{L}_{1} = [10,4,1], \quad \mathcal{L}_{2} = [11,0,4], \quad \mathcal{L}_{3} = [12,9,0], \quad \mathcal{L}_{4} = [13,9,1], \quad \mathcal{L}_{5} = [14,0,7],$$
$$\mathcal{L}_{6} = [15,4,6], \quad \mathcal{L}_{7} = [17,7,7], \quad \mathcal{L}_{8} = [18,6,9], \quad\mathcal{L}_{9} = [22,0,19].$$
Now we need to check whether the above combinatorics can be realized over the real numbers. To this end, first observe that $\mathcal{L}_{1}, ..., \mathcal{L}_{9}$ satisfy the assumptions of Theorem \ref{Schnur}. Combining Shnurnikov`s inequality with $t_{2} = d + t_{5}$ we obtain
\begin{equation}
\label{final}
d - 8 \geq \frac{1}{2}t_{4} + \frac{3}{2}t_{5},
\end{equation}
and it is easy to check that none of $\mathcal{L}_{i}$ satisfies (\ref{final}).
This contradiction finishes the proof.
\end{proof}

Next, we show Theorem B.

\begin{proof}
Suppose that there exists a line configuration $\mathcal{L}$ such that $Y_{4}^{\mathcal{L}}$ is a ball-quotient.
This implies that $\mathcal{L}$ satisfies the following equality:
\begin{equation}
\label{equality11}
9t_{2} + 7t_{3} + t_{4} = 9d + \sum_{r\geq 5} (6r-25)t_{r}.
\end{equation}
Let us recall that Hirzebruch in \cite[p.~140]{Hirzebruch} pointed out that one can improve (\ref{weakHirz}), namely
\begin{equation}
\label{inequality11}
t_{2} + \frac{3}{4}t_{3} \geq d + \sum_{r \geq 5}(r-4)t_{r}.
\end{equation}
Now let us rewrite (\ref{inequality11}) as follows
\begin{equation}
\label{equality12}
9t_{2} + \frac{27}{4}t_{3} \geq 9d + \sum_{r\geq 5}(9r - 36)t_{r}.
\end{equation}
On the other hand, we have
\begin{equation}
\label{inequality12}
9t_{2} + \frac{27}{4}t_{3}  = -t_{4} - \frac{1}{4}t_{3} + 9d + \sum_{r\geq 5} (6r-25)t_{r}.
\end{equation}
Combining (\ref{equality12}) with (\ref{inequality12}) we obtain
\begin{equation}
-t_{4} - \frac{1}{4}t_{3} + 9d + \sum_{r\geq 5} (6r-25)t_{r} \geq 9d + \sum_{r\geq 5}(9r - 36)t_{r},
\end{equation}
which implies $t_{r} = 0$ for $r\geq 3$ and (\ref{equality11}) has the following form
$$t_{2} = d.$$
However, using the combinatorial equality one gets
$$d(d-1) = 2t_{2} = 2d,$$
which implies that either $d = 3$ or $d=0$, a contradiction.
\end{proof}
\begin{remark}
Using almost the same proof one can show that there does not exist any line configuration $\mathcal{L}$ of $d\geq 4$ lines with $t_{d} = t_{d-1} = 0$ such that $Y_{7}^{\mathcal{L}}$ is a ball-quotient.
\end{remark}
Finally, we show Theorem C.

\begin{proof}
Again, our problem boils down to classifying all real line configurations that satisfy the following equality:
\begin{equation}
\label{equality2}
4t_{2} +3t_{3} + t_{4} = 4d + \sum_{r\geq 5}(2r-9)t_{r}.
\end{equation}
It is easy to see that one can automatically exclude the case $t_{d-2} =1$, thus from now on we assume that $t_{d}=t_{d-1} = t_{d-2} = 0$.
Rewriting (\ref{equality2}) in a slightly different way we get
$$t_{2} + \frac{3}{4}t_{3} = d - \frac{1}{4}t_{4} + \sum_{r \geq 5}\left(\frac{1}{2}r - \frac{9}{4}\right)t_{r}.$$
Now combining this with (\ref{strongHirz}), we obtain
$$d - \frac{1}{4}t_{4} + \sum_{r \geq 5}\left(\frac{1}{2}r - \frac{9}{4}\right)t_{r} \geq d + \sum_{r\geq 5}(2r-9)t_{r}$$ and finally
$$-\frac{1}{4}t_{4} \geq \sum_{r\geq 5} \left(\frac{3}{2}r - \frac{27}{4}\right)t_{r}.$$
This implies $t_{r} = 0$ for $r\geq 4$ and it leads to
\begin{equation}
\label{equality3}
t_{2} + \frac{3}{4}t_{3} = d.
\end{equation}
Using the combinatorial equality with (\ref{equality3}) one gets
\begin{equation}
\label{equality4}
\frac{2}{9}d(d-3) = t_{3}.
\end{equation}
On the other hand, by Melchior's inequality
$$t_{2} \geq 3$$
and
$$d(d-1) = 2t_{2} + 6t_{3} \geq 6(1+t_{3}).$$
Now using (\ref{equality4}) we obtain 
$$d^{2} - 9d + 18 \leq 0,$$
which means $d \in \{4,5,6\}$. It is easy to verify now that
all these constrains lead to $d=6, t_{2}=3$ and $t_{3}=4$, which completes the proof.
\end{proof}

\section{Topological $(n_{k})$-configurations}
A topological $(n_k)$ point-line configuration, or simply a topological $(n_k)$-configuration, is a set of $n$ points and $n$ pseudolines in the real projective plane, such that each point is incident with $k$ pseudolines and each pseudoline is incident with $k$ points.  
Much work has been done \cite{Gruenbaum} to study the existence of $(n_k)$-configurations in which all pseudolines are straight lines. In these cases it is useful to know whether there exists at least a topological $(n_k)$-configuration. For $k=4$ the existence of topological $(n_4)$-configurations is known for all $n\geq 17$, see \cite{BokowskiGrunebaumSchewe}. 

Using the inequality of Shnurnikov (\ref{Shnurnikov}), we obtain lower bounds for smallest topological $(n_k)$-configurations for $k>4$.
The corresponding bound for $k=4$ is not sharp and leads to $n \ge 16$, however for $k=5$ not much is known so far.

Now we prove Theorem D.

\begin{proof} When we have a topological $(n_k)$-configuration, we can change the configuration locally (if neccessary) such that $t_s = 0$ for $2 < s < k$ and for $k < s$. This implies that the number of single crossings is

$$t_2 = {n \choose 2} - n\cdot{k \choose 2}$$ 
and the inequality of Shnurnikov becomes

$$ n\cdot (n-1) - n\cdot k \cdot (k-1) > 16 + n\cdot (4\cdot k - 15)$$

$$ n \cdot (n-1 - k \cdot (k-1) - 4 \cdot k + 15) > 16$$ 

$$ n \cdot (n + 14 - k \cdot (k+3)) > 16$$

This implies especially that there are no topological $(n_5)$-configurations for $n < 27$ 
and there are no topological $(n_6)$-configurations for $n < 41$.
\end{proof}

The smallest known topological $(n_5)$-configuration with $n=36$ is due to Leah Wrenn Berman, constructed from two
$(18_4)$-configurations, \cite{Berman}. It will be published elsewhere. An open problem remains to find topological $(n_5)$-configurations for $27 \leq n \leq 35$.

\section*{Acknowledgements}
The second author would like to express his gratitude to Alex K\"uronya, Stefan Tohaneanu and Giancarlo Urz\'ua for very useful conversations on the topic of this paper. Both authors would like to thank Leah Wrenn Berman for her useful suggestions. The project was conducted when the second author was a fellow of SFB $45$ \textit{Periods, moduli spaces and arithmetic of algebraic varieties}, and he was partially supported by National Science Centre Poland Grant 2014/15/N/ST1/02102.

%*****************************************************************************

%***************************************************************************** % Addresses
\bigskip
   Piotr Pokora,
   Instytut Matematyki,
   Pedagogical University of Cracow,
   Podchor\c a\.zych 2,
   PL-30-084 Krak\'ow, Poland.

Current Address:
    Institut f\"ur Algebraische Geometrie,
    Leibniz Universit\"at Hannover,
    Welfengarten 1,
    D-30167 Hannover, Germany. \\
\nopagebreak
   \textit{E-mail address:} \texttt{piotrpkr@gmail.com, pokora@math.uni-hannover.de}
   
\bigskip
J\"urgen Bokowski,
Department of Mathematics, Technische Universit\"at Darmstadt,
Schlossgartenstrasse~7, D-64289 Darmstadt, Germany.
\nopagebreak
   \textit{E-mail address:} \texttt{juergen.bokowski@gmail.com}

%*****************************************************************************

\end{document}